\documentclass[11pt]{article}
\usepackage[iso-8859-7]{inputenc}
\usepackage{a4}
\usepackage{amsfonts}
\usepackage{amsmath}
\usepackage{amsthm}
\usepackage{amssymb,pstricks}

\def\f{\varphi}
\def\n{\nu}
\def\x{\xi}
\def\m{\mu}
\def\l{\lambda}
\def\z{\zeta}
\def\<{\langle}
\def\>{\rangle}
\def\c{\psi}
\def\Z{\mathbb Z}
\def\N{\mathbb N}
\def\ar{\rightarrow}
\def\S{\Sigma}
\def\SS{\mathcal{S}}

\begin{document}

\newtheorem{theorem}{Theorem}[section]
\newtheorem{proposition}[theorem]{Proposition}
\newtheorem{lemma}[theorem]{Lemma}
\newtheorem{corollary}[theorem]{Corollary}
\newtheorem{definition}[theorem]{Definition}
\newtheorem{remark}[theorem]{Remark}
\newtheorem{remar}{Remark}
\newtheorem*{rem}{Remark}
\newtheorem{remarks}[theorem]{Remarks}
\newtheorem*{rems}{Remarks}
\newtheorem{examples}[theorem]{Examples}
\newtheorem{example}[theorem]{Example}
\newtheorem{applications}[theorem]{Applications}
\newcommand{\pr}{\noindent {\it Proof.\quad}}

\title{On the linearity of HNN-extensions with abelian base}
\author{V. Metaftsis, E. Raptis and D. Varsos}
\date{}

\maketitle

\begin{abstract}
We show that an HNN-extension with finitely generated abelian base group is $\Z$-linear if
and only if it is residually finite.  
\end{abstract}

\section{Introduction}

A group $G$ is {\it linear} if it admits a faithful representation into
some matrix group GL$_n(k)$ for some commutative ring $k$. The linearity of groups seems 
to be a difficult property to recognize and in most
cases there is no certain method to do so. The first paper that systematically studied
linearity of groups was Mal'cev \cite{malcev}. Since then many authors have shown the
linearity of certain families of groups but the list of linear groups remains short.
We must mention here that in 1988, Lubotzky (see \cite{lubotzky}) gave necessary and sufficient conditions
for a group to be linear over ${\mathbb C}$. Unfortunately, proving that certain groups satisfy the 
Lubotzky's criterion appears to be not an easy job.

In case we choose $k$ to be the ring of integers $\Z$, the range of examples 
of linear groups shortens even further. The purpose of the present note is to 
investigate the linearity of HNN-extensions of the form $G=\< t, K \mid t a t^{-1}= \f(a), \ a\in A \>$ with $K$ a finitely generated abelian group and $\f: A\ar B$ an isomorphism
between subgroups of $K$. Our main result (Corollary \ref{corollary 1}) shows that those groups
are $\Z$-linear if and only if are residually finite.

Moreover, an interesting side result shows that when a certain isolated subgroup of $K$ is trivial
and $A\cap B\neq 1$ then there is a finite index subgroup $K'$ of $K$ such that the isomorphism 
between $A$ and $B$ is induced by an automorphism of $K'$ of finite order. This allows us
to embed a certain finite index subgroup of the HNN-extension
into a larger group which in turn has a finite index subgroup which is a 
right-angled Artin group. That is enough to prove the linearity of the original
HNN-extension. At the end of the paper we also give various consequences of our main
result concerning HNN-extensions.

We conjecture that fundamental groups of trees of groups with finitely generated 
abelian vertex groups are always linear and even that fundamental groups of graphs with finitely 
generated abelian vertex groups are linear if and only if they are residually finite, but we
have not been able to prove such a result with the techniques at hand.

\section{On the isomorphisms of subgroups of finitely generated free abelian
groups}

Let $G$ be a group and $H$ a subgroup of $G$. The subgroup $H$ is {\it isolated}
in $G$, if whenever $g^n\in H$ for $g\in G$ and $n>0$, then $g\in H$. By $i_G(H)$
we denote the {\it isolated closure} of $H$ in $G$, that is, the intersection of all
isolated subgroups of $G$ that contain $H$. For more on isolated subgroups and the
isolated closure of a group, the reader should consult \cite{b}

For the sequel, $G$ is always the HNN-extension $G=\< t, K\mid tat^{-1}=\f(a),\ a\in A\>$
where $K$ is a finitely generated free abelian group and $A, B$ isomorphic
subgroups of $K$ with $\f : A \ar  B$ the isomorphism induced
by $\f$. Let $D$ be the subgroup of $G$ with $$D = \{ x\in K \mid \mbox{\ for each\ } \nu \in
\mathbb{Z}  \mbox{\ there exists\ } \lambda  = \lambda (\nu) \in
\mathbb{N} \mbox{\ such that\ }  t^{-\nu}x^{\lambda}t^{\nu}\in
K \}.$$ Then $D$ is an isolated subgroup of $K$ (see Proposition 2 in \cite{arv}) and
therefore a direct factor of $K$. Moreover $D  \leq i_{K}(A\cap
B)$ and $i_{K}(A\cap B)/D$ is a free abelian group. In fact $D$ plays a central
r\^ole in the proofs of the main results of \cite{arv,arv2} and apparently
in the present work as well.

We can also describe $D$ as follows: let $M_{0} = A\cap B$,
$M_{1} =\f^{-1}(M_{0})\cap M_{0}\cap \f(M_{0})$ and
inductively $M_{i+1} =\f^{-1}(M_{i})\cap M_{i}\cap
\f(M_{i})$. Then $M_{i+1} \leq M_{i}$ and since $K$ is
finitely generated, there is $k\in \N$ such that 
rank$(M_{k+1}) =$ rank$(M_{k})$. Consequently, 
$D = i_{K}(M_{k})$ (see again \cite{arv,arv2}).  Notice also that if $H = \bigcap_{i = 0}^{\infty}M_{i}$, 
then $H$ contains every subgroup  $L$ of $K$ with $\varphi(L) = L$, and in that 
sense is the largest normal subgroup of $G$ contained in $K$; in other words, $H$
is the core $H=K_G$ of $K$ in $G$.

Finally, we can give an alternative description of $D$ using standard Bass-Serre theory (see \cite{serre}). 
Let $T$ be the standard tree on which $G$ acts. Then $D$ is the subgroup of $K$ such that for every
finite subtree $T'$ of $T$, there is a positive integer $n\in\N$ such that $D^n$ stabilizes
$T'$ pointwise. Notice that the subgroup $H$ defined above is the subgroup of $K$ that
stabilizes the entire tree $T$ pointwise.

\bigskip

In this section we show that if $D=1$ and $A\cap B\neq 1$, then there is an algorithm that 
allows us to consider a certain finite index subgroup $\overline{K}$ of $K$, such that the automorphism 
induced on $\overline{K}$ by $\f$ has finite order.

Fix a generating set for $K$ and hence for $A$ and $B$ and take  a generating element
$a_1$ of $A$ and the powers 
$$\ldots,\f^{-2}(a_1), \f^{-1}(a_{1}),a_{1},\f(a_{1}),\f^2(a_1),\ldots .$$
Since $D = 1$ there are $\lambda_{1},   \mu_{1} \in
\mathbb{N}$ with the property
$$\{   \varphi^{-\lambda _{1}+1}(a_{1}),
\ldots, \varphi ^{-1}(a_{1}), a_{1}, \f(a_{1}), \ldots,
\varphi ^{\mu_{1}-1}(a_{1})\} \subseteq A\cap B,$$
but $\varphi^{-\lambda_{1}}(a_{1}) \in A\backslash (A\cap B)$ and $\varphi ^{\mu _{1}}(a_{1}) \in
B\backslash (A\cap B)$.

Here we can assume that $\langle  \varphi
^{-\lambda _{1}}(a_{1}) \rangle \cap (A\cap B) = 1$ and
$\langle \varphi ^{\mu _{1}}(a_{1}) \rangle \cap (A\cap
B) = 1$. For if $(\varphi ^{-\lambda _{1}}(a_{1}))^{\kappa}\in
A\cap B$ for some $\kappa\in\N$ or  $(\varphi ^{\mu_{1}}(a_{1}))^{\nu}\in A\cap B$ 
for some $\nu\in\N$, then we can replace $a_1$ by $a_1'=(a_{1}
)^{\kappa\cdot\nu} $ and $K$ by a finite index subgroup $K_1$
such that $a_1'$ is a generator of $K_1$ and such that we have (possibly
some other) $\lambda _{1},  \mu _{1} \in \mathbb{N}$ with the
property that they are the greatest positive integers such that
$$\{  \varphi ^{-\lambda _{1}+1}(a_{1}'),  \dots  \varphi
^{-1}(a_{1}'), a_{1}', \f(a_{1}'), \dots , \varphi ^{\mu
_{1}-1}(a_{1}') \} \subseteq A\cap B.$$ Now this procedure
has to stop in finitely many steps. Indeed, if for example
$\< \f^{\m}(a_1)\>\cap (A\cap B)\neq 1$ for all $\m\ge 0$
then $\cap_{k=0}^{\infty}\f^k(A\cap B)\neq 1$ and so
by repeatedly applying $\f^{-1}$ we have that
$\cap_{k=-\infty}^{\infty}\f^{k}(A\cap B)\neq 1$ which
implies that $D\neq 1$, a contradiction to our assumption.

So we can assume that there are maximal $\l_1,\m_1$ with 
$\langle \varphi
^{-\lambda _{1}}(a_{1}) \rangle \cap (A\cap B) = 1$ or
$\langle \varphi ^{\mu _{1}}(a_{1}) \rangle \cap (A\cap
B) = 1$.

The elements of the set
 \begin{equation}\label{chain}
\SS_1=\{  \varphi ^{-\lambda_{1}}(a_{1}), 
\dots  \varphi ^{-1}(a_{1}), a_{1}, \varphi(a_{1}), \dots , \varphi ^{\mu _{1}}(a_{1}) \} \end{equation}
are ( $\mathbb{Z}$-) linearly independent in $K$. Indeed, let 
$$(\f^{-\lambda_1}(a_1))^{\x_{-\l_1}}\ldots(\f^{-1}(a_1))^{\x_{-1}} a_1^{\x_0}
(\f(a_1))^{\x_1}\ldots(\f^{\m_1}(a_1))^{\x_{\m_1}}=1$$ with $\x_j\in\Z$.
 If $j_0$ with 
$-\lambda_{1} \leq j_{0} \leq \mu_{1}  $ is the smallest subscript
such that $\xi _{j_{0}} \not = 0$ we have that
$$\f^{\x_{j_0}j_0}(a_1)=\f^{\x_{j_0+1}j_0+1}(a_1)\ldots\f^{\x_{-1}(-1)}(a_1)
a_1^{\x_0}\f^{\x_1}(a_1)\ldots\f^{\x_{m_1-1}\m_1}(a_1)$$
and from this (by applying $\varphi ^{-\lambda _{1}-j_{0}}$) we have
that 
$$\f^{-\l_{1}\x_{j_0}}(a_1) =\f^{-\l_1-j_0}\Bigl( 
\varphi^{(j_{0}+1)\x_{j_0+1}}(a_1)  \ldots \varphi^{-1\x_{-1}}(a_1)a^{\x_0}_{1} \varphi^{\x_1}(a_1)\dots\f^{\m_1\x_{\m_1}}(a_1)\Bigr)$$ is an element of $A\cap B$, absurd by the 
 fact that $\langle \varphi ^{-\l_{1}}(a_{1}) \rangle \cap (A\cap
B) = 1$.

\bigskip

Now define $\S_1$ to be the subgroup of $K$ generated by the 
linearly  independent elements,
$\S_{1} = \langle \f^{-\l_1}(a_1)=a,\ldots, a_1, \ldots, b=\f^{\m}(a_1) \rangle$
   of rank $\lambda _{1}  + 1 + \mu _{1}$. 

We take $a_2$ to be a generating element of $A$ and repeat the above procedure.
At the end we get a subgroup $\S_2$ of $K$ of finite rank. We continue until we
exhaust all generating elements of $A$ and then we repeat the procedure with
the generating elements of $B$ and $\f^{-1}$. After finitely many steps, we end
up with finitely many sets $\SS_i$, $i=1,\ldots, m$ of 
the form (\ref{chain}) and finitely many corresponding subgroups of $K$, 
$\S_1,\S_2,\ldots,\S_m$.  In fact, we can think of $\SS_i$ as chains of elements
produced by $\f$ and the generators of $A$ and $B$. So view $\SS_1$ as:
$$a\mapsto \f(a) \mapsto \ldots \mapsto \f^{\l_1}(a)=a_1
\mapsto \ldots \mapsto \f^{\xi_1}(a)=b$$ where $\xi_1=\kappa_1+\lambda_1$.

By construction, is obvious that $\< A,B\>=\< \S_1,\S_2,\ldots,\S_m\>$. Moreover, in the
collection $\S_1,\S_2,\ldots,\S_m$ we can define a  partial ordering in the following
way: 
$\S_i \prec \S_j$ if and only if $\S_j$ contains either $\S_i$ or a finite index
subgroup of it. In the first case we delete chain $\SS_i$ from our collection. 
In the second case, we choose finite index subgroup of $K$, such that
$\S_j$ contains $\S_i$ and then we delete $\SS_i$ from our collection. Since we have
finitely many subgroups, after finitely many step, we can find appropriate finite index 
subgroup of $K$ and chains $\SS_i$, $i=1,\ldots, p$ such that $\S_i \nprec\S_j$
for any $i,j\in\{ 1,\ldots, p\}$.

From now on we restrict attention to these remaning chains $\SS_i$, $i=1,\ldots, p$
and their corresponding subgroups $\S_i$.  
By the maximality choice of $\SS_i$ we have that $\S_i\cap\S_j$ with $i\neq j$
is a subgroup of $A\cap B$. 
We choose any $\SS_i$ and extend $\f$ in the following way: for notational simplicity assume 
that $\SS_1$ is still in our collection; define $\f_1$ to be the extension of $\f$ that
   maps $b\mapsto a$. In other words, $\f_1|_A=\f$ (i.e. $\f_1$ is 
   $\f$ compatible) and $\f_1(b)=a$. So $\f_1$ is naturally defined in $\<A,b\>$ with 
   $\f_1(wb^k)=\f(w)a^k$, with $w\in A$. One can easily see that $\f_1$ is
   well defined and is an isomorphism $\f_1: \< A, b\> \ar \< B, a\>$. Moreover, when restricted to $\S_1$ is an automorphism of finite order. Notice that $\f_1$ is now possibly defined automatically on some elements of
$B\setminus (A\cap B)$ other than $b$.

Now use $\f_1$ to re-calculate the remaining chains $\SS_i$. This means that since
we have $\f_1$ defined on $b$ and $\f_1^{-1}$ defined on $a$ we may be able to
extend $\SS_i$ further. For if the last element of the chain $\SS_i$ belongs to $\< A,b\>$ then
$\f_1$ can be applied so to extend the chain further. Similarly, we can apply 
$\f_1^{-1}$ to the first element of $\SS_i$, if that belongs to $\< B,a\>$. 
So let $\SS_i'$ be the chains that are produced after the application of
$\f_1$. 

We claim that none of the chains $\SS_i'$ contains an infinite number
of elements. Indeed, assume 
that $\SS_2$ is one of the remaining chains with 
$$\SS_2=\{x=\f^{-\l_2}(a_2),\ldots, \f^{-1}(a_2),a_2,\f(a_2),\ldots, \f^{\m_2}(a_2)=y\}$$
$$=\{ x, \f(x), \f^2(x), \ldots, \f^{\xi_2}(x)=y\}$$ with $x\in A\setminus (A\cap B)$
and $y\in B\setminus (A\cap B)$ where $\xi_2=\lambda_2+\mu_2$.
 
In order to be able to apply $\f_1$ to $y$ we must have $y\in\< A,b\>$, and since
$y\in B\setminus (A\cap B)$, the only possibility for $y$ is to be of the form
$y=wb^k$ for some $k\in\Z\setminus\{ 0\}$ and some $w\in A\cap B$. 
Now if $w\in\S_1$ then $\S_2\prec\S_1$ a contradiction. If on the other
hand $w\in\S_2$ then $b^k\in\S_2$ and so $\S_1\prec\S_2$, a further
contradiction. Now if $w\in \<\S_1,\S_2\>$ then there are $w_1\in\S_1\setminus \S_2$ and $w_2\in\S_2\setminus\S_1$ such that $w=w_1w_2$. But then $w_1b^k\in \S_1\cap\S_2$ and $w_1b^k\in B\setminus (A\cap B)$, a contradiction to the fact
that $\S_1\cap\S_2$ is a subgroup of $A\cap B$.

Moreover, if $w$ is a word in some other $\S_i$ such that $\f^k(w)=w_1b_1^r$ 
with $w_1\in A\cap B$ then by construction, $w_1'b^r$ is the right end of
$\SS_i$ for some $w_1'\in\S_i$ and the same arguments as before show that 
$w_1\not\in\< \S_1,\S_2,\S_i\>$. Since we have finitely many chains,
 So after finitely many applications of $\f$, $\f^k(w)\not\in\<A,b\>$.  A similar 
argument can be 
applied to $x$. This proves the claim.

\bigskip

We proceed in the same way. That is, for the new chains $\SS_i'$, we calculate the subgroups
generated by their elements and we delete the chains for which their subgroups
are contained in the subgroups of other chains, choosing appropriate finite
index subgroup of $K$ if necessary. We choose a chain randomly and
extend $\f_1$ to $\f_2$ by mapping the last element of this chain to its first element.
Again, $\f_2$ is well defined, $\f_2|_A=\f$ and $\f_2|_{\< A,b\>}$. This procedure terminates in finitely many 
steps and the final $\f_k$ that is defined as an extension of $\f$, is obviously of finite
order. So $\f_k$ is an automorphism of a finite index subgroup of $\< A,B\>$ of 
finite order. We extend $\f_k$
to a finite index subgroup of $K$ by mapping the remaining generators of $K$ to themselves.

So we have shown the following.

\begin{proposition}\label{proposition 1}
Let $G=\< t, K\mid tat^{-1}=\f(a), \ a\in A>$ where $K$ is a finitely generated abelian
base group and $\f:A\ar B$ is the isomorphism induced by $t$. 
Suppose that $D$ is defined as above and that $D = 1.$ Then there exists a finite index subgroup 
$\bar{K}$ of $K$ and a finite order automorphism $\bar{\f}$ of $\bar{K}$, such that  
$\bar{\varphi}\mid _{\bar{K}\cap A} = \varphi$. Moreover, there is an algorithm
which constructs such a $\bar{\f}$ in finitely many steps.  \qed
\end{proposition}

\section{Linearity of HNN-extensions}

\begin{theorem}
 \quad Let $K$ be a finitely generated free abelian group, $\varphi$ an
automorphism of $K$ of finite order and $A$ a subgroup of $K$. Then, the
multiple HNN-extension
\[ G = \langle  t_{1}, \dots , t_{n}, K \mid t_{i}at_{i}^{-1} = 
\varphi (a),  a\in A, i = 1,  \dots , n   \rangle \]
is $\mathbb{Z}$-linear.

\end{theorem}

\pr
We take the HNN-extension
$\tilde{G}$ generated by the elements $\xi_{1},\ldots,\xi_{n}, \zeta, K$, satisfying the following relations $\xi_{i}k\xi_{i}^{-1} = 
\varphi(k),  k\in K, i = 1,  \dots , n,$ $[\zeta, k]=1$ for all $k\in K$.
Notice that each $\x_i$ acts on $K$ as an automorphism of finite order, the same order for
all $i=1,\ldots, n$.

The map  $f : G \longrightarrow \tilde{G}$ with $f(k) = k$
for every element of $K$ and $f(t_{i}) = \xi_{i}\zeta$ defines a
monomorphism. Indeed, the relations in $G$ are preserved by $f$,
so it is a homomorphism, and if $1 \not = g\in G$, then by a
$t$-length argument we can see that $f(g) \not = 1$, namely $f$
is 1-1.

Let $\n$ be the order of the automorphism $\varphi$. Since $\< \x_i, i=1,\ldots, n\>$
generate a free group, we can consider the epimorphism $\c : \< \x_1,\ldots, \x_n\> \ar \Z_{\n}$
with $\c(\x_i)=1$ for all $i=1,\ldots, n$. This epimorphism, extends to an epimorphism of $\tilde{G}$, which
for simplicity we also denote $\c$, $\c: \tilde{G} \ar \Z_{\n}$ by sending all elements of 
$K$ and $\z$ to zero. Let $H$ be the kernel of $\c$. Obviously, $H$ is a subgroup of finite
index in $\tilde{G}$. In order to find a presentation of $H$ we choose a Schreier transversal $U$
for $H$ to be the set $U=\{ 1, \x_1, \x_1^2, \ldots, \x_1^{\n-1}\}$. Then a basis
for $H$ consists of the non-trivial elements of $H$ of the form $u_1\x_i u_2^{-1}$
$i=1,\ldots, n$, $u_1\z u_2^{-1}$ and $u_1 k_ju_2^{-1}$ with $u_1,u_2\in U$
and $k_j$ the generators of $K$. Hence, $H$ is generated by the set
$$\{ \x_1^r\x_i\x_1^{-(r+1)}, \x_1^s\z\x_1^{-s},
 \x_1^sk_j\x_1^{-s},    i=2,\ldots,n, r=1,\ldots, \nu-2, s=0,\ldots, \nu-1,$$
$$\x_1^{\n-1}\x_i, i=1,\ldots, \n-1 \}$$
where $\{ k_j\}$ is a generating set of $K$. Rename the above generating set as follows:
$$x_{ir}=\x_1^r\x_i\x_1^{-(r+1)}$$
$$x_{i\n}=\x_1^{\n-1}\x_i, i=2,\ldots, n, \ \ x=\x_1^{\n}$$
$$z_s=\x_1^s\zeta\x_1^{-s}$$
$$k_{sj}=\x_1^sk_j\x_1^{-s}$$

From the Reidemeister-Schreier rewriting process, the relations of $H$ are 
the relations $uru^{-1}$
where $u\in U$ and $r\in R$, where $R$ is the set of relations of $\tilde{G}$, rewritten 
in terms of the generators of $H$.  So the set of relations for $H$ consists of the following
relations:
$$\x_1^s [\z, k_j] \x_1^{-s}  =[z_s, k_{sj}] \ \ s=1,\ldots, \n-1$$
$$\x_1^s [k_{j_1},k_{j_2}]\x_1^{-s}=[k_{sj_1},k_{sj_2}] \ \ s=1,\ldots, \n-1$$
$$\x_1^s \x_i k_j\x_i^{-1}\x_1^{-s}=\x_1^s\f(k_j)\x_1^{-s}$$
The last set of relations, for $i=1$ becomes
\begin{eqnarray}
k_{(s+1)j}=\f(k_{sj}) & s=0,\ldots, \nu-2\label{relations1}\\
x k_{0j} x^{-1} =\f(k_{(\n-1)j}) & \label{relations2}
\end{eqnarray}
For every $i=2,\ldots,n$ becomes
$$\x_1^s\x_i\x_1^{-(s+1)} \cdot \x_1^{(s+1)} k_j \x_1^{-(s+1)} 
\cdot\x_1^{(s+1)}\x_i^{-1}\x_1^{-s}= \x_1^s \f(k_j) \x_1^{-s}$$ and by renaming
$$x_{is}k_{(s+1)j}x_{is}^{-1}=\f(k_{sj}) \ \ s=0,\ldots, \n-2$$
and
$$x_{i\n} k_{0j}x_{i\n}^{-1}=\f(k_{(\n-1)j}).$$
Now use the set of relations (\ref{relations1}) and (\ref{relations2}) to
solve for $k_{sj}$ and use a Tietze transormation to replace $k_{sj}$ to the above relations to get a new set of relations for our group
$$k_{(s+1)j}=\f(k_{sj}) \ \   s=0,\ldots, \nu-2$$
$$k_{0j} =x^{-1} \f(k_{(\n-1)j})x$$
$$[\f(k_{sj_1}),\f(k_{sj_2})] \ \ s=0,\ldots, \n-1$$
$$[z_s, \f(k_{(s-1)j})] \ \  s=1,\ldots, \n-1$$
$$[z_0, x^{-1}\f(k_{(\n-1)j})x]$$
$$[x_{is},\f(k_{sj})],  \ \ s=0,\ldots, \n-2$$
$$x_{i\n}x^{-1} \f(k_{(\n-1)j})xx_{i\n}^{-1} =  \f(k_{(\n-1)j}), \ \  i=1,\ldots, n$$
Since $\f$ is an automorphism, every $\f(k_{sj})$ is a generating set of a copy of $K$
for every $s=0,\ldots,\n-1$ and so we can replace the generating set $k_{sj}$
by $\f(k_{sj})$ and eliminate $k_{sj}$. Moreover, we can use Tietze trasfromations
to replace $x_{i\nu}$ by $x_{i\n}x^{-1}=y_i$, $i=1,\ldots,n$ and $z_0$ by
$xz_ox^{-1}=z_0'$. Then the a set of relations for the kernel $H$ is
$$[\f(k_{sj_1}),\f(k_{sj_2})] \ \ s=0,\ldots, \n-1$$
$$[z_s, \f(k_{(s-1)j})] \ \  s=1,\ldots, \n-1$$
$$[z_0', \f(k_{(\n-1)j})]$$
$$[x_{is},\f(k_{sj})],  \ \ s=0,\ldots, \n-2$$
$$[y_i, \f(k_{(\n-1)j})], \ \  i=1,\ldots, n$$
But the above implies that $H$ has a presentation where all relations are
commutators of the generators and therefore is a right-angled Artin group. 
Consequently $H$ is $\Z$-linear. (For more on right-angled 
Artin groups and its linearity the reader can see \cite{humphries,hw,dj}.)  But it
is known that the linearity is closed under taking  finite extensions
or subgroups. Therefore the group $\tilde{G}$ is $\mathbb{Z}$-linear
and hence, so if $G$. \qed

The above technique can actually show that if $K$ is a right-angled Artin group with standard generating
set $S$, and  $S_1, S_2$ are subsets of $S$ then any HNN-extension $\< t, S \mid tS_1t^{-1}=S_2\>$ is 
$\Z$-linear. This is only a special case of a more general result by Hsu and Leary \cite{hl}. 
\vspace{.3in}

\begin{proposition}
 Let $K$ be a finitely generated free abelian group and
$A, B$ isomorphic subgroups of $K$ with $\varphi :
A \longrightarrow  B$ an isomorphism. Suppose that $D = 1$,
then the $HNN$-extension
\[ G = \langle  t,  K \mid tat^{-1} = \f(a), \ a\in A\>\]
is $\mathbb{Z}$-linear.
\end{proposition}

\pr
 From the proposition \ref{proposition 1}  there exists a finite index
subgroup $\bar{K}$ in $K$ and an automorphism $\bar{\varphi}$ of
$\bar{K}$ of finite order such that $\bar{\varphi}\mid
_{\bar{K}\cap A} = \varphi $. Let
$G_{1} = \langle t, \bar{K}  \rangle ^{G}$ be the normal
closure of $\langle t, \bar{K}  \rangle$ in $K$. Evidently
$G_{1}$ is of finite index in $G$. Let
$k_{1}, k_{2}, \dots , k_{n}$ be representatives of $\bar{K}$
in $K$. We take $t_{i} = k_{i}tk_{i}^{-1},
  i = 1, 2, \dots , n$, then using the Schreier rewriting
process we obtain for $G_{1}$ the presentation
$G_{1} = \langle 
t_{1}, \dots , t_{n}, \bar{K} \mid t_{i}at_{i}^{-1} = 
\bar{\varphi} (a),  a\in \bar{K}\cap A, i = 1,  \dots , n 
\rangle$. The group $G_{1}$ is $\mathbb{Z}$-linear by the previous
theorem. So $G$ is, as finite extension of $G_{1}$. \qed

\begin{proposition} \label{proposition 2}
 Let $K$ be a finitely generated abelian group and $A, B$ isomorphic subgroups
of $K$ with $\varphi : A \longrightarrow  B$ an isomorphism.
Suppose that $D$ is finite, then the $HNN$-extension
\[ G = \langle  t,  K \mid tat^{-1} = 
\varphi (a),\  a\in A    \rangle \]
is $\mathbb{Z}$-linear.
\end{proposition}

\pr
 The group $G$ is residually finite, since $D$ is finite (see \cite{arv}). Therefore there 
exists a finite index normal subgroup $N$
such that $N\cap D = 1$. Let $K_{1} = K\cap N$,
  $A_{1} = A\cap N,   B_{1} = B\cap N$ and $\varphi _{1}
:A_{1} \longrightarrow B_{1}$ the induced isomorphism (the group
$N$ is normal in $G$). It is clear that the corresponding subgroup
$D_{1} = \{ x\in K_{1} \mid \mbox{\ for each\ } \nu \in
\mathbb{Z}  \mbox{\ there exists\ } \lambda  = \lambda (\nu) \in
\mathbb{N} \mbox{\ such that\ }  t^{-\nu}x^{\lambda}t^{\nu}\in
K_{1} \}$ is trivial. So, the HNN-extension $\<t, K_1\mid tat^{-1}=\f_1(a),\ a\in A_1\>$
satisfies the hypotheses of Proposition
\ref{proposition 1}, consequently there is a finite
index subgroup $K_{2}$ of $K_{1}$ and an automorphism
$\bar{\varphi} _{1}$ of $K_{2}$ of finite order which extends
$\varphi _{1}$. The normal closure $\langle t, K_{2} 
\rangle^{G}$ is of finite index in $G$ and, as in the previous
proposition, we get that $\langle t, K_{2}  \rangle^{G}$ (and therefore
the group $G$) is $\mathbb{Z}$-linear. \qed

Let $K$ be a finitely generated abelian group and $A, B$ isomorphic subgroups
of $K$ with $\varphi : A \longrightarrow  B$ an isomorphism and
\[ G = \langle  t,  K \mid tat^{-1} = 
\varphi (a),  a\in A    \rangle \] be the corresponding
HNN-extension. Let also $H$ be the largest subgroup of $K$ such that $\varphi
(H) = H$, i.e. the core of $K$ in $G$.

\begin{theorem} \label{theorem 1}

Let $K$ be a f.g. abelian group and $A, B$ proper, isomorphic 
subgroups of $K$ with $\varphi : A \longrightarrow  B$ an
isomorphism and 
\[ G = \langle  t,  K \mid tat^{-1} = 
\varphi (a),  a\in A    \rangle \] be the corresponding
$HNN$-extension.
The group $G$ is $\mathbb{Z}$-linear if and only if the group
$\bar{G} = G/H$ is $\mathbb{Z}$-linear.
\end{theorem}

\pr
Assume $G$ to be linear. Then $G$ is residually finite, therefore,
by the main Theorem in \cite{arv}, we have that $H$ is of finite index in
$D$. The quotient $G/H$ has the HNN-presentation
$G/H = \langle t, K \mid t (A/H)
t^{-1} = B/H,  \varphi_{H}  \rangle$, where $\varphi _{H}$ is
the induced isomorphism (since $\varphi (H) = H$). The
corresponding subgroup $D_{H}$ is equal to $\bar{D} = D/H$ which
is finite. Therefore by the previous proposition $G/H$ is linear.

For the converse, since the group $\bar{G} = G/H$ is linear,
there is an homomorphism $\vartheta : G \longrightarrow R$ to a
linear group $R$ with Ker$\vartheta  \leq H$. On the other hand
the linearity of $\bar{G} = G/H$ implies the residually
finiteness of it, but the corresponding
$\bar{H} = (K/H)_{\bar{G}}$ is trivial. Therefore, by the main
Theorem in \cite{arv} the group $\bar{D} = D/H$ must be finite, which,
by the  Theorem 2 in \cite{arv2}, implies that there exists a finitely generated abelian
group $X$ such that $K\leq X$ and an automorphism $\bar{\varphi}$
of $X$ with $\bar{\varphi}|_{A} = \varphi$.
Now the obvious homomorphism $\varrho :
G \longrightarrow X \rtimes \langle \bar{\varphi}  \rangle$ is an
embedding on $K$, so Ker$\varrho \cap K = 1$. The linearity of
$G$ follows from the linearity of the groups $R$,
$X \rtimes \< \bar{\varphi}  \>$ (which is linear since it is polycyclic \cite{segal}) 
and the fact that
Ker$\vartheta \cap$ Ker$\varrho = 1$. \qed

\begin{corollary} \label{corollary 1}
Let $K$ be a f.g. abelian group and $A, B$ proper isomorphic 
subgroups of $K$ with $\varphi : A \longrightarrow  B$ an
isomorphism and
\[ G = \langle  t,  K \mid tat^{-1} = 
\varphi (a),  a\in A  \rangle \] be the corresponding
$HNN$-extension.
The group $G$ is $\mathbb{Z}$-linear if and only if it is
residually finite.
\end{corollary}

\pr Suppose that $G$ is residually finite. By the main Theorem in \cite{arv}
we have that $H$ is of finite index in $D$, where $H$ and $D$ are
the subgroups of $K$ defined above. Therefore the group $G/H$ is
linear by proposition \ref{proposition 2}. The result now follows
from the previous theorem.

For the converse, it is well known that a finitely generated linear group is
residually finite. \qed

\section{Some further thoughts}
\begin{remar}
 In the proof of Theorem \ref{theorem 1}, in order to prove that the
linearity of $G$ implies the linearity of $G/H$ it is worth
remarking that we used only that $G$ is residually finite and
proposition \ref{proposition 2}, which in turn depends on
proposition \ref{proposition 1}. We can give a direct proof using
heavily that $G$ is linear as follows.
\end{remar}

The quotient $G/H$ has the HNN-presentation
$G/H = \langle t, K \mid t A/H
t^{-1} = B/H,  \varphi_{H}  \rangle$, where $\varphi _{H}$ is
the induced isomorphism (since $\varphi (H) = H$). The
corresponding subgroup $D_{H}$ is equal to $\bar{D} = D/H$ which
is finite, so $G/H$ is residually finite. Therefore there exists a
normal subgroup $N$ of finite index in $G$ such that $H\leq N$ and
$N\cap D = H$. By the structure theorem of Bass-Serre theory (see \cite{serre}), 
the structure of the subgroups of an HNN-extension $N$ is the fundamental group of a finite graph
of groups with vertex groups of the form $N\cap gKg^{-1}$,   $g\in
G$ and edge groups of the form $N\cap gAg^{-1}$,   $g\in G$. Since
$H = K_{G}$ and $N$ is normal in $G$, $H$ is contained in every
edge (and vertex) group. Now the group $D$ is isolated in $K$, which 
implies that $D\cap N$ is isolated in $K\cap N$. So
$H = gHg^{-1}$ is isolated (therefore a direct factor)  in every
vertex group. This means, by the normal form of elements of $N$,
that $N = H\ltimes M$ for some subgroup $M$ of the linear group $G$.
Consequently $N/H \simeq M$ is linear. But $N/H$ is of finite
index in $G/H$, so $G/H$ is linear.

\begin{remar}
 The \emph{value} of the propositions \ref{proposition 1} and
\ref{proposition 2} consists in the fact that exhibiting  an
internal property of finitely generated abelian groups we obtain a
\emph{tangible} criterion for the linearity of HNN-extensions with
base group a finitely generated abelian group.
\end{remar}

\begin{corollary}
Let $K$ be a finitely generated abelian group, $A,  B$ proper
subgroups of $K$ and $\varphi :A \longrightarrow B$ an
isomorphism.

\begin{enumerate}
 \item The subgroup $D$ is finite.
\item There exists a finite index subgroup $K_{1}$ of $K$ and
an automorphism $\varphi _{1}$ of $K_{1}$ of finite order which
extends $\varphi$ in the sense that $\varphi _{1_{\mid A\cap
K_{1}}} = \varphi$.
\item There exists an abelian group $X$ which contains the
group $K$ as a subgroup of finite index and an automorphism
$\vartheta$ of $X$ of finite order such that $\vartheta _{\mid
A} = \varphi$.
\end{enumerate}
For the above statements we have the following:

1) implies 2).

2) is equivalent to 3).

1) implies 3).
 \end{corollary}

\pr 
1)\quad implies \quad 2). This is the proposition
\ref{proposition 2}.

2) is equivalent to 3). Assuming 2), it is easy to 
construct the group $X$ adding roots to the group $K_{1}$.

Assuming 3), we take $K_{1} = X^{n}$, where $n$ is the index
of $K$ in $X$.

1)\quad implies \quad 3) is a consequence of the above two. \qed

\begin{remar}
In account of Theorem 2 in \cite{arv2}, we see that there is a (weak)
equivalence where there is no needed for the automorphism $\vartheta$ to
have finite order. On the other hand, the statements 2) or
3) do not imply 1). For example the free abelian group
$K = \langle a,  b \mid [ a, b ] \rangle$ with
$A = \langle a \rangle = B$ and $\varphi$ the identity
satisfies (trivially) both 2) and 3), but not 1).
\end{remar}

In Corollary \ref{corollary 1} it is assumed that the associated
subgroups $A$ and $B$ are proper subgroups in the base group $K$.
In the case where one of them is all the base group (e.g.
$K = A$), then the HNN-extension
$G = \langle t, K \mid tat^{-1} = \f(a), \ a\in K  \rangle$ is a
residually finite group, as a solvable constructible group, it is
$\mathbb{Q}$-linear. (This result and a good account of basic
properties of solvable constructible groups can be found in
\cite{strebel}). Then $G$ is not a $\mathbb{Z}$-linear group. This is concluded
form the fact that the subgroup $B$ is not closed in the profinite
topology of $G$ (see e.g. in \cite{rv}), on the other hand if $G$ was
$\mathbb{Z}$-linear, then by Theorem 5 p. 61 in \cite{rtv}, all subgroups
of $K$ must be closed in the profinite topology of $G$.

\bigskip

In the case where the base group of an HNN-extension is not a f.g.
abelian group only miscellaneous cases are known for the linearity
of these groups.

\begin{proposition}
Let $K$ be any finitely generated linear group and $\varphi$ an
isomorphism between finite subgroups $A$ and $B$ of $K$. The
HNN-extension $G = \langle t, K \mid tat^{-1} = \varphi
(a),  a\in A  \rangle$ is linear.
\end{proposition}

\pr At first, $K$ is residually finite and since $A$ and $B$ are
finite groups, the group $G$ is residually finite. Therefore there exists
a normal subgroup $N$ of finite index in $G$ such that $N\cap
A = N\cap B = 1$. From the structure theorem of Bass-Serre theory, 
the subgroup $N$ is a free product of a free
group and a finite family of subgroups of kind $N\cap gKg^{-1}$,
  $g\in G$. Then the result of Nisnevi\v c (\cite{nisnevich}, see also \cite{wehrfritz}) 
implies that if the group $K$ is linear of degree $d\geq 2$, then
$N$ is linear of degree at most $d+1$. So $G$ is linear. \qed

\begin{proposition}
(cf. Theorem 1.3. in \cite{rtv}) The HNN-extension
$G = \langle t, K \mid tat^{-1} = \varphi (a),  a\in A 
\rangle$ with base group $K$ a polycyclic-by-finite group and
proper associated subgroups $A$ and $B = \varphi (A)$ of finite
index in $K$ is $\mathbb{Z}$-linear if and only if it is subgroup
separable.
\end{proposition}

In the case where the associated subgroups are not of finite index
we have a simple (very) special result.

\begin{proposition}
Let $K$ be the split extension of a polycyclic-by-finite group $A$
by a polycyclic-by-finite group $C$ ($K = A\rtimes C$) and $\f:A\ar A$ an
automorphism of $A$. Then the
HNN-extension
$G = \langle t, K \mid tat^{-1} = \f(a), \ a\in A\rangle$
is linear.
\end{proposition}

\pr Evidently $G = A\rtimes\langle t, C  \rangle$. Therefore
$G/A \simeq \< t  \> \ast C$, so it is linear. On
the other hand the map $\vartheta : G \longrightarrow 
A\rtimes\langle \varphi  \rangle$ which sends every element $a\in A$
to itself, every element of $C$ to the trivial element  and $t$ to
$\varphi$ is a well defined homomorphism with $ker \vartheta \cap
A = 1$. Therefore, $G$ is linear. \qed

\vspace{.3in}

Department of Mathematics, University of the Aegean, Karlovasi 832 00, Samos, Greece. {\it Email: vmet@aegean.gr}.

Department of Mathematics, University of Athens, Panepistimiopolis 157 84, Athens, Greece. 

{\it Email of E. Raptis: eraptis@math.uoa.gr}. 

{\it Email of D. Varsos: dvarsos@math.uoa.gr}.

\end{document}